\documentclass[11pt]{article}
\usepackage{amsfonts}
\usepackage{mathrsfs}

\usepackage[latin1]{inputenc}
\usepackage{amsmath,amssymb}
\usepackage{latexsym}
\usepackage[active]{srcltx}
\usepackage{multicol}  
\usepackage{multirow} 
\usepackage{threeparttable}  
\usepackage[
bookmarks=true,         
bookmarksnumbered=true, 
colorlinks=true, pdfstartview=FitV, linkcolor=blue, citecolor=blue,
urlcolor=blue]{hyperref}

\newtheorem{theorem}{Theorem}[section]

\newtheorem{lemma}[theorem]{Lemma}

\newtheorem{remark}[theorem]{Remark}
\numberwithin{equation}{section}
\def\qed{{\hfill $\square$ \bigskip}}
\def\R{{\mathbb R}}
\def\E{{{\mathbb E}\,}}

\def\N{{\mathbb N}}

\def\Var{{\mathop {{\rm Var\, }}}}

\def\square{{\vcenter{\vbox{\hrule height.3pt
        \hbox{\vrule width.3pt height5pt \kern5pt
           \vrule width.3pt}
        \hrule height.3pt}}}}
\def \eref#1{\hbox{(\ref{#1})}}

\begin{document}

\title{Kernel entropy estimation for long memory linear processes with infinite variance}
\date{\today}
\author{Hui Liu and Fangjun Xu\thanks{F. Xu is partially supported by National Natural Science Foundation of China (Grant No.11871219, No.11871220).} \\
}
\maketitle

\begin{abstract}  

Let $X=\{X_n: n\in\mathbb{N}\}$ be a long memory linear process with innovations in the domain of attraction of an $\alpha$-stable law $(0<\alpha<2)$. Assume that the linear process $X$ has a bounded probability density function $f(x)$. Then, under certain conditions, we consider the estimation of the quadratic functional $\int_{\mathbb{R}} f^2(x) \,dx$ by using the kernel estimator 
\[
T_n(h_n)=\frac{2}{n(n-1)h_n}\sum_{1\leq j<i\leq n}K\left(\frac{X_i-X_j}{h_n}\right).
\] 
The simulation study for long memory linear processes with symmetric $\alpha$-stable innovations is also given.

\noindent


\end{abstract}

\section{Introduction}

Let $X=\{X_n:\, n\in\mathbb{N}\}$ be a linear process defined by 
\begin{equation}\label{def-Xn}
X_n=\sum_{i=0}^{\infty}a_{i}\varepsilon_{n-i},
\end{equation}
where the innovations $\varepsilon_i$ are i.i.d. real-valued random variables belonging to the domain of attraction of an $\alpha$-stable law ($0<\alpha<2$), the coefficients  $a_0=1$, $a_i\sim c_0i^{-\beta}$, $i=1,2,\cdot\cdot\cdot,$ and $c_0$ is a positive constant. Here $a_i\sim a'_i$ means that $a_i/a'_i\to1$ as $i\to\infty$. By Kolmogorov three-series theorem, the linear process $X$ in \eref{def-Xn}  converges almost surely if $\alpha\beta>1$. Assume that the linear process $X$ has a bounded probability density function $f(x)$. Then the study of quadratic functional $\int_{\R} f^2(x)\,dx$ will help us to get more information on entropies related to the linear process $X$, say, quadratic R\'{e}nyi entropy $R(f)=-\ln (\int_{\R} f^2(x)\,dx)$ and Shannon entropy $S(f)=-\int_{\R} f(x)\ln f(x)\,dx$.

Entropy is widely applied in the fields of information theory, statistical classification, pattern recognition and so on since it is a measure of uncertainty in a probability distribution. In the literature, different estimators for the quadratic functional and entropies with independent data  have been well studied. However, there are very few works on estimations of the quadratic functional and the corresponding entropies for dependent case.  In \cite{KLS},  K\"{a}llberg,  Leonenko and Seleznjev extended the $U$-statistics method to $m$-dependence sequence. They showed the rate optimality and asymptotic normality of the $U$-statistics estimator for multivariate sequence.  In \cite{A}, Ahmad obtained the strong consistency of the quadratic functional by orthogonal series method for stationary time series with strong mixing condition. In \cite{Sang-Xu}, kernel entropy estimation for quadratic functional and related entropies of regular time series data under certain mild conditions were studied. Although the linear processes in \cite{Sang-Xu} can have infinite variance, it can only deal with the short memory case, see Definition 2.1 and Example 3.2 in \cite{Sang-Xu} for more details. To the best of our knowledge, general results for quadratic functional estimations and related entropies of long memory linear processes with infinite variance are still unknown.

In this paper, for the linear process $X=\{X_n:\, n\in\mathbb{N}\}$ defined in (\ref{def-Xn}), we only focus on the case $1<\alpha \beta<2$.  According to Definition 2.1 in \cite{Sang-Xu}, this corresponds to the long memory case. When innovations are symmetric $\alpha$-stable random variables, one can also refer to \cite{Hsing-1999} for the definition of such long memory linear processes. To estimate the quadratic functional $\int_{\R} f^2(x)\, dx$ of the linear process $X=\{X_n:\, n\in\mathbb{N}\}$ defined in $\eref{def-Xn}$, we shall apply the kernel method
\begin{align} \label{def-Tn}
T_n(h_n)=\frac{2}{n(n-1)h_n} \sum_{1\le j< i\le n}K\left(\frac{X_i-X_j}{h_n}\right),
\end{align}
where the kernel $K$ is a symmetric and bounded function with $\int_{\R} K(u)\, du = 1$ and $\int_{\R} u^2|K(u)|\, du<\infty$. The bandwidth sequence $h_n$ satisfies $0<h_n\to 0$ as $n\to\infty$.

Throughout this paper, if not mentioned otherwise, the letter $c$ with or without a subscript, denotes a generic positive finite constant whose exact value is independent of $n$ and may change from line to line. We use $\iota$ to denote the imaginary unit $\sqrt{-1}$. For a complex number $z$, we use $\overline{z}$ and $|z|$ to denote its conjugate and modulus, respectively. For any integrable function $g(x)$, its Fourier transform is defined as~$\widehat{g}(u)=\int_{\mathbb{R}}e^{\iota x u}g(x)\, dx$. Moreover, we let $\phi(\lambda)$ be the characteristic function of linear process $X=\{X_n:\, n\in\mathbb{N}\}$,  and $\phi_{\varepsilon}(\lambda)$ the characteristic function of innovations. That is, $\phi(\lambda)=\mathbb{E}[e^{\iota \lambda X_n}]$ and $\phi_{\varepsilon}(\lambda)=\mathbb{E}[e^{\iota \lambda \varepsilon_{1}}]$.  For simplicity of notation, we always assume that the coefficients $a_i$ in the definition of the linear process $X$ are nonzero.

The paper has the following structure. The main results are given in Section 2.  A simulation study is given in Section 3. Section 4 is devoted to the proofs of Theorems \ref{thm1} and \ref{thm2} based on the Fourier transform and the projection method.

\bigskip

\section{Main Results}

It is well known that the characteristic function of an $\alpha$-stable law $S_{\alpha}(\sigma,\eta,\mu)$ has the form
\[
e^{\iota\lambda \mu -\sigma^{\alpha}|\lambda|^{\alpha}(1-\iota \eta \operatorname{sign}(\lambda) \omega(\lambda,\alpha))},
\]
where $0<\alpha \leq 2,\sigma>0,-1 \leqslant \eta \leqslant 1,\mu \in \mathbb{R}$, and
\begin{equation*}\label{def-omega}
\omega(\lambda,\alpha)= \begin{cases} \tan(\frac{\pi \alpha}{2})
	 & \text { for } \alpha \neq 1 \\ 
	 2 \pi^{-1} \log|\lambda|
	 & \text { for } \alpha=1.\end{cases}
\end{equation*}
It is called symmetric (or $S\alpha S$) if $\eta=\mu=0$, and standard if $\sigma=1$. For more details on stable laws, we refer to \cite{Samorodnitsky and Taqqu-1994}.  Let $\xrightarrow{\mathcal{L}}$ denote the convergence in distribution. A random variable $Y$ is said to be in the domain of attraction of an $\alpha$-stable law if there exist i.i.d. random variables $Y_n$ with the same distribution as $Y$, real numbers $A_n$, and strictly positive numbers $B_n$ such that 
\[
\frac{\sum\limits^n_{i=1}Y_i-A_n}{B_n}\xrightarrow{\mathcal{L}} S_{\alpha}(\sigma,\eta,\mu)
\] 
as $n\to \infty$, see, for example, \cite{Ibragimov and Linnik-1971}.  Since the innovations $\varepsilon_i$  belong to the domain of attraction of an $\alpha$-stable law, by Theorem 2.6.5 in \cite{Ibragimov and Linnik-1971}, the characteristic function $\phi_{\varepsilon}$ of $\varepsilon_1$ satisfies 
\[
|\phi_{\varepsilon}(\lambda)|=e^{-c_{\alpha} |\lambda|^{\alpha} L(\lambda)(1+o(1))}
\]
as $\lambda\to 0$, where $c_{\alpha}$ is a positive constant depending on $\alpha$ and $L(\lambda)$ is a slowly varying function as $\lambda\to 0$. Since the coefficients satisfy $a_0=1$, $a_i \sim c_0i^{-\beta}$, it is easy to see that 
\[
\sum\limits_{i=0}^{\infty}\sqrt{\Var(e^{\iota \lambda a_i\varepsilon_1})}=\infty\quad \text{but}\quad \sum\limits_{i=0}^{\infty}\Var(e^{\iota \lambda a_i\varepsilon_1})<\infty
\] 
for $1<\alpha\beta<2$. So, by Definition 2.1 in \cite{Sang-Xu},  the linear proces $X=\{X_n: n\in\mathbb{N}\}$ defined in \eref{def-Xn} has long memory. This is consistent with the definition of long memory linear processes when innovations are symmetric $\alpha$-stable random variables, see \cite{Hsing-1999}.

Let $G(x)$ be the distribution function of $\varepsilon_1$. We introduce the following assumptions on the innovation $\varepsilon_1$.
\begin{enumerate}
\item[({\bf A1})] There exist strictly positive constants $c_1$ and $\delta$ such that $|\phi_{\varepsilon}(\lambda)|\leq \frac{c_1}{1+|\lambda |^{\delta}}$ for all  $\lambda\in\mathbb{R}$;

\item[({\bf A2})] There exists non-negative constants $c_{-}$ and $c_+$ with $c_{-}+c_+>0$ such that 
\begin{equation*} 
	\lim_{x\to-\infty}|x|^{\alpha}G(x)=c_{-} \quad\text{and}\quad \lim_{x\to+\infty}x^{\alpha}(1-G(x))=c_+;
\end{equation*}

\item[({\bf A3})] $G(x)$ is twice differentiable with derivatives $G^{(j)}(x) (j=1,2)$ satisfying the following inequalities: for any $x, y \in \mathbb{R},|x-y| \leqslant 1$, $j=1,2$
\[
|G^{(j)}(x)| \leq C(1+|x|)^{-\alpha}
\]
and
\[
|G^{(j)}(x)-G^{(j)}(y)| \leq C|x-y|(1+|x|)^{-\alpha} .
\]
\end{enumerate}

The assumption ({\bf A1}) and $a_i\sim c_0 i^{-\beta}$ imply that  (i) the linear process $X$ defined in \eref{def-Xn} has a bounded probability density function $f(x)$; (ii) the characteristic function $\phi(\lambda)$ of $X_n$ decays at any polynomial rate to $0$ as $|\lambda|\to+\infty$. Moreover, ({\bf A1}) implies that there exists $m\in\mathbb{N}$ such that $|\phi_{\varepsilon}(\lambda)|^m\leq \frac{c_2}{1+|\lambda |^4}$. Following the proof of Lemma 1 of \cite{Giraitis et al-1996} we can get that $f(x)$ is twice continuously differentiable and all its derivatives up to the second order are uniformly bounded, see also the {\bf P1} in \cite{Honda-2009b}. 

Assumptions ({\bf A2}) and ({\bf A3}) are only needed in Theorem \ref{thm2} to obtain the desired limiting theorems. The assumption ({\bf A2}) shows that $\varepsilon_1$ belongs to the domain of attraction of $\alpha$-stable law. For more details on the domain of attraction of $\alpha$-stable law, we refer to \cite{Ibragimov and Linnik-1971}. If $\varepsilon_1$ is a symmetric $\alpha$-stable random variable, then its distribution function $G$ satisfies all assumptions ({\bf A1}), ({\bf A2}) and ({\bf A3}).

The following are main results of our paper.

\begin{theorem}\label{thm1}
Assume that  ({\bf A1}) holds. Then, for any $\eta\in(0,\alpha-\frac{1}{\beta})$, there exist positive constants $c_1$ and $c_2$ depending on $\eta$ such that
\begin{align*}\label{case1-mean}
\Big|\E T_n(h_n)-\int_{\mathbb{R}}f^2(x)\, dx\Big|\leq c_1\left(n^{1-(\alpha-\eta)\beta}+h_n^2\right)
\end{align*}
and
\begin{align*}
\E\Big|T_n(h_n)-\E T_n(h_n)-\frac{1}{n}\sum^n_{i=1}Y_i\Big|\leq c_2\left(n^{1-(\alpha-\eta)\beta}+\frac{1}{\sqrt{n^3 h^2_n}}+\frac{1}{\sqrt{n^2 h_n}}+h^2_nn^{\frac{1-(\alpha-\eta)\beta}{2}} \right),
\end{align*}
where $Y_i=2\Big( f(X_i)-\int_{\mathbb{R}} f^2(x)\, dx \Big)$.
\end{theorem}

\begin{theorem}\label{thm2}
Under the assumptions of Theorem \ref{thm1} and $nh_n\to \infty$ as $n\to\infty$,
\begin{enumerate}
\item[(1)] if ({\bf A2}) and ({\bf A3}) hold with $c_-=c_+$, $1<\alpha<2$, $\dfrac{1}{\alpha}<\beta<1$ and $\lim\limits_{n\to\infty}n^{\frac{(\alpha\beta-1)(2-\alpha)}{4\alpha}+\frac{\eta\beta}{4}}h_n=0$ for $\eta\in(0,\frac{(\alpha\beta-1)(\alpha-1)}{\alpha\beta})$, then we have 
\begin{equation}\label{asym1}
	n^{\beta-\frac{1}{\alpha}}\left[T_n(h_n)-\E T_n(h_n)\right]\xrightarrow{\mathcal{L}} \tilde{c} Z,
\end{equation}
where $Z$ is a standard $S\alpha S$ random variable and  
\[
\tilde{c}=2c_{0}\left(2 c_+ \frac{\Gamma(2-\alpha) \cos (\alpha \pi / 2)}{1-\alpha} \int_{-\infty}^{1} \int_{0}^{1}(t-s)_{+}^{-\beta} dt ds\right)^{1 / \alpha}\int_{\mathbb{R}}f(x)df(x) .
\]

\item[(2)] if ({\bf A2}) and ({\bf A3}) hold with $c_-=c_+$, $1<\alpha<2$, $1<\beta<\dfrac{2}{\alpha}$ and $\lim\limits_{n\to\infty}n^{\frac{(\alpha\beta-1)(2-\alpha\beta)}{4\alpha\beta}+\frac{\eta\beta}{4}}h_n=0$ for $\eta\in(0,\frac{(\alpha\beta-1)^2}{\alpha\beta^2})$, then we have 
	\begin{equation}\label{asym2}
		n^{1-\frac{1}{\alpha\beta}}\left[T_n(h_n)-\E T_n(h_n)\right]\xrightarrow{\mathcal{L}} c_+^{1/\alpha\beta}c_{f}^{+}L^{+}+c_-^{1/\alpha\beta} c_{f}^{-}L^{-},
	\end{equation}
	where $L^{-}$and $L^{+}$ are i.i.d. random variables with stable law $S_{\alpha \beta}(1,1,0)$ and
	\begin{align*}
			c_{f}^{\pm}=2 \tilde{\sigma} \int_{0}^{\infty}\left(f_{\infty}(\pm u)-f_{\infty}(0)\right) u^{-(1+1 / \beta)} \mathrm{d}u
	\end{align*}
with $f_{\infty}(x)=\E[f\left(X_{1}+x\right)]$ and $\tilde{\sigma}=\left\{\frac{c_{0}^{\alpha}(\alpha \beta-1)}{\Gamma(2-\alpha \beta)|\cos (\pi \alpha \beta / 2)| \beta^{\alpha \beta}}\right\}^{1 /(\alpha \beta)}$.

\item[(3)] if ({\bf A2}) holds, $0<\alpha<1$, $1<\alpha\beta<2$ and $\lim\limits_{n\to\infty}n^{\frac{(\alpha\beta-1)(2-\alpha\beta)}{4\alpha\beta}+\frac{\eta\beta}{4}}h_n=0$ for $\eta\in(0,\frac{(\alpha\beta-1)^2}{\alpha\beta^2})$, then we have 
	\begin{equation}\label{asym3}
		n^{1-\frac{1}{\alpha\beta}}\left[T_n(h_n)-\E T_n(h_n)\right]\xrightarrow{\mathcal{L}}c_+^{1/\alpha\beta} c_{f}^{+}L^{+}+c_-^{1/\alpha\beta}c_{f}^{-}L^{-},
	\end{equation}
	where $c_f^{\pm}$ and $L^{\pm}$ are defined in (2).
\end{enumerate}	
\end{theorem}

\begin{remark} 
The bandwidth $h_n$ in Theorem \ref{thm2} can be chosen independent of $\alpha$ and $\beta$. Note that for any $x\in(1,2)$, $g(x)=\frac{(x-1)(2-x)}{4x}\in (0, \frac{3-2\sqrt{2}}{4}]$ and $\frac{3-2\sqrt{2}}{4}<\frac{3-2.8}{4}=\frac{1}{20}$.  Choose $\eta\in(0,\alpha-\frac{1}{\beta})$ small enough. We only need to require $n^{\frac{1}{20}}h_n\to 0$ and $nh_n\to\infty$ as $n\to\infty$. Moreover, $\E T_n(h_n)$ in $\eref{asym1}$, $\eref{asym2}$ and $\eref{asym3}$ can be replaced by $\int_{\mathbb{R}}f^2(x)\, dx$ if $h_n$ satisfies $\lim\limits_{n\to\infty}n^{\frac{\alpha\beta-1}{2\alpha}}h_n=0$ in $\eref{asym1}$  and $\lim\limits_{n\to\infty}n^{\frac{\alpha\beta-1}{2\alpha\beta}}h_n=0$ in $\eref{asym2}$ and $\eref{asym3}$. Clearly, if $h_n=O(n^{-\frac{1}{4}})$ and $nh_n\to \infty$ as $n\to\infty$, we can replace $\E T_n(h_n)$ in Theorem $\ref{thm2}$ by $\int_{\mathbb{R}}f^2(x)\, dx$.
\end{remark}

\bigskip

\section{Simulation}

We carry out a simulation study to examine properties of the kernel entropy estimator for the linear process $X=\{X_n:\, n\in\mathbb{N}\}$ defined in $\eref{def-Xn}$. Here we assume that the innovation $\varepsilon_1$ follows the standard symmetric $\alpha$-stable law $S_{\alpha}(1,0,0)$. Moreover, we take the usual normal kernel function $K(x)=\frac{1}{\sqrt{2\pi}}e^{-\frac{x^2}{2}}$, the bandwidth $h_n=n^{-\frac{1}{5}}$ and coefficients 
\[
a_i= \begin{cases} 1,
		& \text { for } i=0\\ 
		 i^{-\beta},
		& \text { for } i\geq 1\end{cases}.
\]
Since the innovations $\varepsilon_i$ are i.i.d. $S_{\alpha}(1,0,0)$ random variables with $0<\alpha<2$, we have $\phi_{\varepsilon}(u)=e^{-|u|^{\alpha}}$.  Therefore, the characteristic function of linear process $X=\{X_n:\, n\in\mathbb{N}\}$ is written as
\[
\phi(u)=\E[e^{\iota\lambda X_n}]=\prod_{i=0}^{\infty}\E[e^{\iota\lambda a_i \varepsilon_{n-i}}]=e^{-|u|^{\alpha}\sum\limits_{i=0}^{\infty}|a_i|^{\alpha}}.
\] 
Now, by Plancherel theorem, the quadratic functional of the linear process  $X=\{X_n:\, n\in\mathbb{N}\}$ can be obtained as
\[
\int_{\mathbb{R}} f^2(x)\, dx=\frac{1}{2\pi}\int_{\mathbb{R}} |\phi(\lambda)|^2\, d\lambda=\dfrac{1}{\pi \alpha \Big(2\sum\limits_{i=0}^{\infty}|a_i|^{\alpha}\Big)^{\frac{1}{\alpha}} }\Gamma(\frac{1}{\alpha}).
\]

We perform the simulation study by using the software MATLAB.  Here, the sample sizes are $n=1000, 2000, 5000$ and we always simulate $N=1000$ times. In the following two tables, Mean, Var and Mse stand for the sample means, the sample variances and the sample mean squared errors, respectively.

The true values of $\int_{\mathbb{R}}f^2(x)\, dx$ and simulation results are summarized in the following two tables. Table 1 is for $\alpha=0.5$, while Table 2 is for $\alpha=1.5$. From these tables, we observe that
\begin{itemize}
\item[(i)]  As $n$ increases,  the estimated value of $\int_{\mathbb{R}} f^2(x)\, dx$ approaches the true value and the bias steadily decreases;     
\item[(ii)] The estimator performs pretty well if $\alpha\beta$ is close to $2$.
\end{itemize}

\begin{table}[htb] \label{table1}
	\centering
	\caption{$\alpha=0.5$}  \medskip
	\begin{tabular}{ccccccclccccc}  
		\hline  
		\multicolumn{3}{c}{\multirow{2}{*}{$\beta$}}&
		\multicolumn{3}{c}{\multirow{2}{*}{True value}}& &
		\multicolumn{1}{c}{$n$}& & 1000&2000&5000&\\
		\multicolumn{3}{c}{}&\multicolumn{3}{c}{}& &\multicolumn{1}{c}{$h_n=n^{-1/5}$}& &0.2513 &0.2187&0.1821&\\   
		\hline  
		\multicolumn{3}{c}{\multirow{3}{*}{2.5}}&\multicolumn{3}{c}{\multirow{3}{*}{0.0051}}& &
		Mean& & 0.0073& 0.0070& 0.0065&\\    
		\multicolumn{3}{c}{}&	\multicolumn{3}{c}{}&  &
		Var($\times 10^{-3}$)& & 0.0099&0.0074&0.0046\\
		\multicolumn{3}{c}{}&	\multicolumn{3}{c}{}& &
		Mse($\times 10^{-3}$)& & 0.0148&0.0111&0.0064&\\
		\hline  
		\multicolumn{3}{c}{\multirow{3}{*}{3.5}}&\multicolumn{3}{c}{\multirow{3}{*}{0.0181}}& &
		Mean& & 0.0184& 0.0185& 0.0182&\\    
		\multicolumn{3}{c}{}&	\multicolumn{3}{c}{}& &
		Var($\times 10^{-3}$)& & 0.0147&0.0080&0.0037&\\
		\multicolumn{3}{c}{}&	\multicolumn{3}{c}{}& &
		Mse($\times 10^{-3}$)& & 0.0148&0.0081&0.0037&\\
		\hline  
		\multicolumn{3}{c}{\multirow{3}{*}{3.9}}&\multicolumn{3}{c}{\multirow{3}{*}{0.0219}}& &
		Mean& & 0.0219& 0.0221& 0.0219&\\    
		\multicolumn{3}{c}{}&	\multicolumn{3}{c}{}& &
		Var($\times 10^{-3}$)& & 0.0146&0.0068&0.0028&\\
		\multicolumn{3}{c}{}&	\multicolumn{3}{c}{}& &
		Mse($\times 10^{-3}$)& & 0.0146&0.0068&0.0028&\\
		\hline  
	\end{tabular}
\end{table}

\begin{table}[htb] \label{table2}
	\centering
	\caption{$\alpha=1.5$}  \medskip
	\begin{tabular}{ccccccclccccc}  
		\hline  
		\multicolumn{3}{c}{\multirow{2}{*}{$\beta$}}&
		\multicolumn{3}{c}{\multirow{2}{*}{True value}}& &
		\multicolumn{1}{c}{$n$}& & 1000&2000&5000&\\
		\multicolumn{3}{c}{}&\multicolumn{3}{c}{}& &\multicolumn{1}{c}{$h_n=n^{-1/5}$}& &0.2513 &0.2187&0.1821&\\   
		\hline  
		\multicolumn{3}{c}{\multirow{3}{*}{0.9}}&\multicolumn{3}{c}{\multirow{3}{*}{0.0668}}& &
		Mean& & 0.0738& 0.0728& 0.0705&\\    
		\multicolumn{3}{c}{}&	\multicolumn{3}{c}{}&  &
		Var($\times 10^{-3}$)& & 0.1110&0.0713&0.0544\\
		\multicolumn{3}{c}{}&	\multicolumn{3}{c}{}& &
		Mse($\times 10^{-3}$)& & 0.1601&0.1076&0.0679&\\
		\hline  
		\multicolumn{3}{c}{\multirow{3}{*}{1.1}}&\multicolumn{3}{c}{\multirow{3}{*}{0.0840}}& &
		Mean& & 0.0858& 0.0857& 0.0846&\\    
		\multicolumn{3}{c}{}&	\multicolumn{3}{c}{}& &
		Var($\times 10^{-3}$)& & 0.0777&0.0420&0.0247&\\
		\multicolumn{3}{c}{}&	\multicolumn{3}{c}{}& &
		Mse($\times 10^{-3}$)& & 0.0807&0.0447&0.0251&\\
		\hline  
		\multicolumn{3}{c}{\multirow{3}{*}{1.3}}&\multicolumn{3}{c}{\multirow{3}{*}{0.0935}}& &
		Mean& & 0.0939& 0.0940& 0.0935&\\    
		\multicolumn{3}{c}{}&	\multicolumn{3}{c}{}& &
		Var($\times 10^{-3}$)& & 0.0525&0.0268&0.0124&\\
		\multicolumn{3}{c}{}&	\multicolumn{3}{c}{}& &
		Mse($\times 10^{-3}$)& & 0.0527&0.0270&0.0124&\\
		\hline  
	\end{tabular}
\end{table}

\bigskip

\section{Proofs}

In this section, we will prove Theorems \ref{thm1} and \ref{thm2}. To begin with, we introduce the projection method and two lemmas. Lemma \ref{lem1} is on the characteristic function of innovations $\varepsilon_i$ and Lemma \ref{lem2} gives the desired estimation of the covariance 
\[
\E\left[ (e^{\iota \lambda X_i}-\phi(\lambda))(e^{-\iota \lambda  X_j}-\phi(-\lambda))\right].
\]
For each $i \in \mathbb{Z}$, let $\mathcal{F}_{i}$ be the $\sigma$-field generated by random variables $\left\{\varepsilon_{k}: k \leq i\right\}$. Given an integrable complex-valued random variable $Z$, we define the projection operator $\mathcal{P}_{i}$ as
\[
\mathcal{P}_{i} Z=\E[Z|\mathcal{F}_{i}]-\E[Z|\mathcal{F}_{i-1}]
\]
for each $i \in \mathbb{Z}$. It is easy to see that $\E[\mathcal{P}_{i} Z \mathcal{P}_{j} W]=0$ if $i\neq j$, $\mathbb{E} |Z|^2<\infty$ and $\mathbb{E} |W|^2<\infty$.

\begin{lemma}  \label{lem1}
If $\varepsilon$ is in the domain of attraction of an $\alpha$-stable law with $\alpha\in (0,2)$ and $\phi_{\varepsilon}(\lambda)$ is its characteristic function, then for any $\eta\in (0,\alpha)$, there exists a positive constant  $c_{\alpha,\eta}$ such that
\[
\E\big|e^{ \iota \lambda\varepsilon}-\phi_{\varepsilon}(\lambda)\big|^2\leq c_{\alpha,\eta}\left(|\lambda|^{\alpha-\eta}\wedge 1\right).
\]
\end{lemma}

\noindent
{\it Proof}:  By Theorem 2.6.5 in \cite{Ibragimov and Linnik-1971}, 
\[
|\phi_{\varepsilon}(\lambda)|^2
=e^{-c_1 |\lambda|^{\alpha}L(\lambda)(1+o(1))}
\] 
as $\lambda\to 0$, where $L$ is a slowly varying function as $\lambda\to 0$. Therefore,
\[
\E|e^{ \iota \lambda \varepsilon_1}-\phi_{\varepsilon}(\lambda)|^2=1-\left|\phi_{\varepsilon}(\lambda)\right|^2\leq c_{\alpha,\eta}\left(|\lambda|^{\alpha-\eta}\wedge 1\right),
\]
where in the last inequality we used the fact that~$|1-e^{-x}|\leq x$~for~$x\geq 0$ and $\lim\limits_{\lambda\to 0}|\lambda|^{\eta}L(\lambda)=0$ for any $\eta>0$. \qed

\begin{lemma}  \label{lem2}
For any $1\leq i\neq j\leq n$ and $\eta\in(0,\alpha-\frac{1}{\beta})$, there exists a positive constant  $c_{\eta}$ such that
\[
\bigg|\E\Big[(e^{\iota \lambda X_i}-\phi(\lambda))(e^{-\iota \lambda  X_j}-\phi(-\lambda))\Big]\bigg|\leq \frac{c_{\eta}}{1+|\lambda|^4} |i-j|^{1-(\alpha-\eta)\beta}.
\]
\end{lemma}

\noindent
{\it Proof}: For any $i\geq 1$,
\begin{align}\label{decomposition}
e^{\iota \lambda X_i}-\phi(\lambda)
&=\sum^{i}_{k=-\infty}\mathcal{P}_{k}(e^{\iota \lambda X_i}-\phi(\lambda)) \nonumber \\
&=\sum^{i}_{k=-\infty}\left(\prod^{i-k-1}_{\ell=0}\phi_{\varepsilon}(\lambda a_{\ell})\right)\left(e^{\iota \lambda a_{i-k} \varepsilon_k}-\phi_{\varepsilon}(\lambda a_{i-k})\right) e^{\iota \lambda\sum\limits^{\infty}_{\ell=i-k+1}a_{\ell}\varepsilon_{i-\ell}}.
\end{align}
It suffices to consider the case $i>j$. 
Using the decomposition \eref{decomposition}, we can obtain that
\begin{align*}
\E\left[ (e^{\iota \lambda X_i}-\phi(\lambda))(e^{-\iota \lambda  X_j}-\phi(-\lambda))\right]=\sum^{j}_{k=-\infty} \text{I}_1\times \text{I}_2 \times \text{I}_3,
\end{align*}
where 
 \begin{align*}
\text{I}_1&=\prod^{i-k-1}\limits_{\ell=0}\phi_{\varepsilon}(\lambda a_{\ell})\prod\limits^{j-k-1}_{\ell=0}\phi_{\varepsilon}(-\lambda a_{\ell}),\\
\text{I}_2&=\mathbb{E}\big[(e^{\iota \lambda a_{i-k} \varepsilon_k}-\phi_{\varepsilon}(\lambda a_{i-k}))(e^{-\iota \lambda a_{j-k} \varepsilon_k}-\phi_{\varepsilon}(-\lambda a_{j-k}))\big],\\
\text{I}_3&=\prod^{\infty}_{\ell=1}\phi_{\varepsilon}(\lambda(a_{i-k+\ell}-a_{j-k+\ell}))
\end{align*}
with the convention $\prod\limits^{-1}_{\ell=0}\phi_{\varepsilon}(\lambda a_{\ell})=1$.

By Assumption ({\bf A1}), there exists $m_0\in\N$ such that $\prod\limits^{m_0}_{\ell=0}\left|\phi_{\varepsilon}(\lambda a_{\ell})\right|$ is less than a constant multiple of $\frac{1}{1+|\lambda |^6}$. Hence, for the case $i-j>m_0$ or $k<j-m_0$,  $|\text{I}_1|$ is less than a constant multiple of $\frac{1}{1+|\lambda |^6}$. Moreover, for the case $1\leq i-j\leq m_0$ and $j-m_0\leq k\leq j$, $\lim\limits_{n\to\infty}n^{\beta}a_n=c_0$ implies that there exist infinitely many $\ell (\geq 1)$ such that 
\[
a_{i-k+\ell}-a_{j-k+\ell}=a_{i-j+j-k+\ell}-a_{j-k+\ell}\neq 0.
\] 
So, by Assumption ({\bf A1}),  $|\text{I}_3|$ is less than a constant multiple of $\frac{1}{1+|\lambda |^6}$ in the case $1\leq i-j\leq m_0$ and $j-m_0\leq k\leq j$. 

Therefore, using Cauchy-Schwartz inequality and Lemma \ref{lem1} with $\eta\in(0, \alpha-\frac{1}{\beta})$,
\begin{align*}
&\Big|\E\left[ (e^{\iota \lambda X_i}-\phi(\lambda))(e^{-\iota \lambda  X_j}-\phi(-\lambda))\right]\Big|\\
&\leq \frac{c_1}{1+|\lambda |^6}\sum^{i\wedge j}_{k=-\infty}\Big|\mathbb{E}\big[(e^{\iota \lambda a_{i-k} \varepsilon_k}-\phi_{\varepsilon}(\lambda a_{i-k}))(e^{-\iota \lambda a_{j-k} \varepsilon_k}-\phi_{\varepsilon}(-\lambda a_{j-k}))\big]\Big|\\
&\leq \frac{c_2}{1+|\lambda |^6}|\lambda|^{\alpha-\eta}\sum^{i\wedge j}_{k=-\infty}|a_{i-k}|^{\frac{\alpha-\eta}{2}}|a_{j-k}|^{\frac{\alpha-\eta}{2}}\\
&\leq \frac{c_3}{1+|\lambda |^4} \sum^{\infty}_{\ell=1}\big(|i-j|+\ell\big)^{-\frac{(\alpha-\eta)\beta}{2}}\ell^{-\frac{(\alpha-\eta)\beta}{2}} \\
&\leq \frac{c_4}{1+|\lambda |^4} \, |i-j|^{1-(\alpha-\eta)\beta}\int^{\infty}_0(1+x)^{-\frac{(\alpha-\eta)\beta}{2}}x^{-\frac{(\alpha-\eta)\beta}{2}}dx \\
&\leq \frac{c_5}{1+|\lambda |^4}\, |i-j|^{1-(\alpha-\eta)\beta}.
\end{align*}
This gives the desired estimate. \qed
 
Now we give the proof of Theorem \ref{thm1}.

\noindent
{\it Proof of Theorem \ref{thm1}} The proof will be done in several steps. 

\noindent
{\bf Step 1.} We give the estimation for $\big|\E T_n(h_n)-\int_{\mathbb{R}}f^2(x)\mathrm{d}x\big|$. Using the Fourier inverse transform, we can obtain that
\begin{align*}
&\E T_n(h_n)-\int_{\R} f^2(x)\, dx\notag\\
&=\frac{1}{\pi n(n-1)}\sum_{1\leq j<i\leq n}\int_{\R}\widehat{K}(\lambda h_n) \, \E\big[ (e^{\iota \lambda X_i}-\phi(\lambda))(e^{-\iota \lambda  X_j}-\phi(-\lambda))\big]\, d\lambda\\
&\qquad\qquad\qquad+\frac{1}{2\pi}\int_{\R}\big(\widehat{K}(\lambda h_n)-\widehat{K}(0)\big)|\phi(\lambda)|^2\, d\lambda\\
&=: \text{II}_1+\text{II}_2.
\end{align*}

By Lemma \ref{lem2} and the boundedness of $\widehat{K}$, we can obtain that
\begin{align*}
|\text{II}_1|\leq \frac{c_1}{n^2} \sum_{1\leq j<i\leq n}|i-j|^{1-(\alpha-\eta)\beta}\leq c_2\, n^{1-(\alpha-\eta)\beta}.
\end{align*}

Since $K$ is a symmetric and bounded function with $\int_{\R} K(u)\, du=1$ and $\int_{\R}u^2 |K(u)|\, du<\infty$,
\begin{align*}
|\widehat{K}(\lambda h_n)-\widehat{K}(0)|=\left|\int_{\R} (e^{\iota\lambda h_n u}-1)K(u) \, du\right|\leq |\lambda|^2 h^2_n \int_{\R} u^2 |K(u)|\, du.
\end{align*} 
Then, by Assumption ({\bf A1}), $|\text{II}_2|$ is less than a constant multiple of $h^2_n$. Hence
\begin{align*}
\Big| \E T_n(h_n)-\int_{\R} f^2(x)\, dx \Big|\leq c_3\left(n^{1-(\alpha-\eta)\beta}+h^{2}_n\right).
\end{align*}

{\bf Step 2.} We give the decomposition for $T_n(h_n)-\E T_n(h_n)$. It is easy to see that
\begin{align} \label{decomp}
T_n(h_n)-\E T_n(h_n)=2N_n+A_{n}-\E A_{n}+B_{n}-\E B_{n},
\end{align}
where
\begin{align*}
N_n&=\frac{1}{\pi n(n-1)}\sum_{1\leq j<i\leq n,\, i-j>m}\int_{\R} \widehat{K}(\lambda h_n) \big(e^{\iota \lambda X_i}-\phi(\lambda) \big) \phi(-\lambda) \, d\lambda,
\end{align*}
\begin{align*}
A_{n}&=\frac{1}{n(n-1)h_n} \sum_{1\le j<i\le n,\, i-j\leq m}K\left(\frac{X_i-X_j}{h_n}\right),
\end{align*}
\begin{align*}
B_{n}&=\frac{1}{\pi n(n-1)}\sum_{1\leq j<i\leq n,\,  i-j>m}\int_{\R} \widehat{K}(\lambda h_n)\big(e^{\iota \lambda  X_i}-\phi(\lambda) \big)\big(e^{-\iota \lambda  X_j}-\phi(-\lambda) \big) \, d\lambda
\end{align*}
and the proper choice of the natural number $m$ will be specified in {\bf Step 4.}

{\bf Step 3.} We estimate $\E |A_{n}-\mathbb{E} A_{n}|^2$. For $1\leq j<i\leq n$, we observe that  
\[
X_i-X_j=\sum^{\infty}_{\ell=0} a^{i,j}_{\ell}\epsilon_{i-\ell},
\] 
where
\begin{equation} \label{aij}
a^{i,j}_{\ell}= \begin{cases} a_{\ell},
		& \text { for } 0\leq\ell\leq i-j-1,\\ 
		 a_{\ell}-a_{\ell-(i-j)},
		& \text { for } \ell\geq i-j.\end{cases}
\end{equation}

For $1\leq j<i\leq n$ with $i-j\leq m$, it is easy to see that  there exists $m_1\in\mathbb{N}$ such that $m_1\geq m$ and 
\begin{align} \label{m1}
\prod\limits^{m_1}_{\ell=0}|\phi_{\varepsilon}(a^{i,j}_{\ell}\lambda)|\leq \frac{c_4}{1+|\lambda|^4}.
\end{align}
Let
\[
A^{i,j}_n=\frac{1}{n(n-1)h_n}\left(K\Big(\frac{X_{i}-X_{j}}{h_n}\Big)-\mathbb{E}K\Big(\frac{X_{i}-X_{j}}{h_n}\Big)\right).
\] 
Then
\begin{align} \label{An}
&\E |A_{n}-\mathbb{E} A_{n}|^2 \nonumber\\
&=\sum_{\substack{ 1\leq j_1<i_1\leq n,\, 1\leq j_2<i_2\leq n;\\ i_1-j_1\leq m, \, i_2-j_2\leq m,\, |i_2-i_1|\leq 2m_1}} \mathbb{E}[A^{i_1, j_1}_nA^{i_2, j_2}_n]+\sum_{\substack{ 1\leq j_1<i_1\leq n, \, 1\leq j_2<i_2\leq n;\\ i_1-j_1\leq m, \, i_2-j_2\leq m, \, |i_2-i_1|>2m_1}} \mathbb{E}[A^{i_1, j_1}_nA^{i_2,j_2}_n].
\end{align}
Boundedness of the kernel function $K$ implies that 
\begin{align}\label{An1}
\sum_{\substack{ 1\leq j_1<i_1\leq n, \, 1\leq j_2<i_2\leq n;\\ i_1-j_1\leq m,\,  i_2-j_2\leq m, \, |i_2-i_1|\leq 2m_1}} \big|\mathbb{E}[A^{i_1,j_1}_nA^{i_2, j_2}_n]\big|\leq \frac{c_5}{n^3h^2_n}.
\end{align}
To estimate
\[
\sum_{\substack{ 1\leq j_1<i_1\leq n, \, 1\leq j_2<i_2\leq n;\\ i_1-j_1\leq m, \, i_2-j_2\leq m, \, |i_2-i_1|>2m_1}} \mathbb{E}[A^{i_1, j_1}_nA^{i_2, j_2}_n],
\]
it suffices to consider the case $i_2-i_1>2m_1$. 

Applying the projection operator $\mathcal{P}_k$ on the terms $A^{i,j}_n$, we see that
\begin{align*} 
\mathbb{E}[A^{i_1,j_1}_nA^{i_2,j_2}_n]
&=\sum^{i_1}_{k=-\infty}\mathbb{E}[\mathcal{P}_kA^{i_1,j_1}_n\mathcal{P}_kA^{i_2,j_2}_n] \nonumber\\
&=\sum^{i_1}_{k=i_1-m_1}\mathbb{E}[\mathcal{P}_kA^{i_1,j_1}_n\mathcal{P}_kA^{i_2,j_2}_n]+\sum^{i_1-m_1-1}_{k=-\infty}\mathbb{E}[\mathcal{P}_kA^{i_1,j_1}_n\mathcal{P}_kA^{i_2,j_2}_n].
\end{align*}

By Fourier inverse transform and the boundedness of $\widehat{K}$,
\begin{align} \label{|Aij|}
&|\mathbb{E}[\mathcal{P}_kA^{i_1,j_1}_n\mathcal{P}_kA^{i_2,j_2}_n]| \nonumber\\
=&\Big|\frac{1}{4\pi^2n^2(n-1)^2}\int_{\mathbb{R}^2} \widehat{K}(\lambda_1h_n)\widehat{K}(\lambda_2h_n)\nonumber\\
&\quad\times \mathbb{E}\big[\mathcal{P}_k (e^{\iota\lambda_1(X_{i_1}-X_{j_1})}-\mathbb{E}e^{\iota\lambda_1(X_{i_1}-X_{j_1})})\mathcal{P}_k (e^{\iota\lambda_2(X_{i_2}-X_{j_2})}-\mathbb{E}e^{\iota\lambda_2(X_{i_2}-X_{j_2})}) \big]\, d\lambda_1 d\lambda_2\Big|\nonumber\\
\leq &\frac{c_6}{n^4}\int_{\mathbb{R}^2} \left|\mathbb{E}\big[\mathcal{P}_k (e^{\iota\lambda_1(X_{i_1}-X_{j_1})}-\mathbb{E}e^{\iota\lambda_1(X_{i_1}-X_{j_1})})\mathcal{P}_k (e^{\iota\lambda_2(X_{i_2}-X_{j_2})}-\mathbb{E}e^{\iota\lambda_2(X_{i_2}-X_{j_2})}) \big]\right|\, d\lambda_1 d\lambda_2.
\end{align}
In the case $i_1-m_1\leq k\leq i_1$, recall the choice of $m_1$ in $\eref{m1}$ and $i_2-i_1>2m_1$, we see that 
\[
\left|\mathbb{E}\big[\mathcal{P}_k (e^{\iota\lambda_1(X_{i_1}-X_{j_1})}-\mathbb{E}e^{\iota\lambda_1(X_{i_1}-X_{j_1})})\mathcal{P}_k (e^{\iota\lambda_2(X_{i_2}-X_{j_2})}-\mathbb{E}e^{\iota\lambda_2(X_{i_2}-X_{j_2})}) \big]\right|
\]
in $\eref{|Aij|}$ is less than a constant multiple of 
\[
\left|\mathbb{E}\big[(e^{\iota\lambda_1a^{i_1,j_1}_{i_1-k}\varepsilon_k}-\phi_{\varepsilon}(\lambda_1a^{i_1,j_1}_{i_1-k}))(e^{\iota\lambda_1a^{i_2,j_2}_{i_2-k}\varepsilon_k}-\phi_{\varepsilon}(\lambda_2a^{i_2,j_2}_{i_2-k}))\big] \right|\frac{1}{1+|\lambda_2|^4}\prod^{\infty}_{\ell=m_1+1}|\phi_{\varepsilon}(a^{i_1,j_1}_{\ell}\lambda_1+a^{i_2,j_2}_{\ell}\lambda_2)|.
\]
According to Assumption ({\bf A1}) and $1\leq i_{\theta}-j_{\theta}\leq m$ for $\theta=1,2$, it is easy to see that
\[
\int_{\mathbb{R}}\prod^{\infty}_{\ell=m_1+1}|\phi_{\varepsilon}(a^{i_1,j_1}_{\ell}\lambda_1+a^{i_2,j_2}_{\ell}\lambda_2)|\, d\lambda_1\leq M,
\]
where $M$ a finite positive number independent of $\lambda_2, i_1, j_1, i_2, j_2$.

In the case $k\leq i_1-m_1-1$, recall the choice of $m_1$ in $\eref{m1}$ and $i_2-i_1>2m_1$, we see that 
\[
\left|\mathbb{E}\big[\mathcal{P}_k (e^{\iota\lambda_1(X_{i_1}-X_{j_1})}-\mathbb{E}e^{\iota\lambda_1(X_{i_1}-X_{j_1})})\mathcal{P}_k (e^{\iota\lambda_2(X_{i_2}-X_{j_2})}-\mathbb{E}e^{\iota\lambda_2(X_{i_2}-X_{j_2})}) \big]\right|
\]
in $\eref{|Aij|}$ is less than a constant multiple of 
\[
\left|\mathbb{E}\big[(e^{\iota\lambda_1a^{i_1,j_1}_{i_1-k}\varepsilon_k}-\phi_{\varepsilon}(\lambda_1a^{i_1,j_1}_{i_1-k}))(e^{\iota\lambda_2a^{i_2,j_2}_{i_2-k}\varepsilon_k}-\phi_{\varepsilon}(\lambda_2a^{i_2,j_2}_{i_2-k}))\big] \right|\frac{1}{1+|\lambda_1|^4}\frac{1}{1+|\lambda_2|^4}.
\]
Therefore, by using Cauchy-Schwartz inequality, Lemma $\ref{lem1}$ and the definition of $a^{i,j}_{\ell}$ in $\eref{aij}$, we can obtain that
\begin{align*}
\sum\limits^{i_1}_{k=i_1-m_1}\left| \mathbb{E}[\mathcal{P}_kA^{i_1,j_1}_n\mathcal{P}_kA^{i_2,j_2}_n]\right|
&\leq \frac{c_7}{n^4}\sum\limits^{i_1}_{k=i_1-m_1} \int_{\mathbb{R}^2}
|a^{i_2,j_2}_{i_2-k}\lambda_2|^{\frac{\alpha-\eta}{2}}\frac{1}{1+|\lambda_2|^4}\prod^{\infty}_{\ell=m_1+1}|\phi_{\varepsilon}(a^{i_1,j_1}_{\ell}\lambda_1+a^{i_2,j_2}_{\ell}\lambda_2)|\, d\lambda_1d\lambda_2\\
&\leq \frac{c_8}{n^4}\sum\limits^{i_1}_{k=i_1-m_1} \int_{\mathbb{R}}
|a^{i_2,j_2}_{i_2-k}\lambda_2|^{\frac{\alpha-\eta}{2}}\frac{1}{1+|\lambda_2|^4}\, d\lambda_2\\
&\leq \frac{c_9}{n^4} \sum\limits^{i_1}_{k=i_1-m_1} |a^{i_2,j_2}_{i_2-k}|^{\frac{\alpha-\eta}{2}}\\
&\leq \frac{c_{10}}{n^4} |i_2-i_1|^{-\frac{(\alpha-\eta)\beta}{2}}
\end{align*}
and 
\begin{align*}
\sum^{i_1-m_1-1}_{k=-\infty}\left|\mathbb{E}[\mathcal{P}_kA^{i_1,j_1}_n\mathcal{P}_kA^{i_2,j_2}_n]\right|
&\leq \frac{c_{11}}{n^4}\sum^{i_1-m_1-1}_{k=-\infty}\int_{\mathbb{R}^2} |a^{i_1,j_1}_{i_1-k}\lambda_1|^{\frac{\alpha-\eta}{2}}|a^{i_2,j_2}_{i_2-k}\lambda_2|^{\frac{\alpha-\eta}{2}}\frac{1}{1+|\lambda_1|^4}\frac{1}{1+|\lambda_2|^4} \, d\lambda_1d\lambda_2\\
&\leq \frac{c_{12}}{n^4}\sum^{i_1-m_1-1}_{k=-\infty}|a^{i_1,j_1}_{i_1-k}|^{\frac{\alpha-\eta}{2}}|a^{i_2,j_2}_{i_2-k}|^{\frac{\alpha-\eta}{2}}\\
&\leq \frac{c_{13}}{n^4}\sum^{i_1-m_1-1}_{k=-\infty}|a_{i_1-k}|^{\frac{\alpha-\eta}{2}}|a_{i_2-k}|^{\frac{\alpha-\eta}{2}}\\
&\leq \frac{c_{14}}{n^4}|i_2-i_1|^{1-(\alpha-\eta)\beta}.
\end{align*}
Hence 
\begin{align} \label{An2}
&\sum_{\substack{ 1\leq j_1<i_1\leq n, \, 1\leq j_2<i_2\leq n;\\ i_1-j_1\leq m,\, i_2-j_2\leq m, \, |i_2-i_1|>2m_1}} \left|\mathbb{E}[A^{i_1,j_1}_nA^{i_2,j_2}_n]\right| \nonumber\\
&\leq \frac{c_{15}}{n^4}\sum_{\substack{ 1\leq j_1<i_1\leq n, \, 1\leq j_2<i_2\leq n;\\ i_1-j_1\leq m,\, i_2-j_2\leq m,\, |i_2-i_1|>2m_1}} (|i_2-i_1|^{-\frac{(\alpha-\eta)\beta}{2}}
+|i_2-i_1|^{1-(\alpha-\eta)\beta})\nonumber\\
&\leq \frac{c_{16}}{n^{1+(\alpha-\eta)\beta}}.
\end{align}
Combining $(\ref{An})$, $(\ref{An1})$ and $(\ref{An2})$ gives  
\begin{align} \label{Ane}
\E |A_{n}-\mathbb{E} A_{n}|^2\leq c_{17}\left(\frac{1}{n^3h^2_n}+\frac{1}{n^{1+(\alpha-\eta)\beta}}\right).
\end{align}

{\bf Step 4.} We estimate $\E|B_{n}-\mathbb{E} B_{n}|$. For each $i\in\mathbb{N}$ and $\lambda\in\mathbb{R}$, define
\[
H(X_i)(\lambda)=e^{\iota \lambda  X_{i}}-\phi(\lambda).
\] 
Then 
\begin{align} \label{B}
B_{n}=B_{n,1}+B_{n,2}+B_{n,3},
\end{align}
where 
\begin{align*}
B_{n,1}
&=\frac{1}{\pi n(n-1)}\sum_{1\leq j<i\leq n,\, i-j>m}\int_{\R} \widehat{K}(\lambda h_n)\sum^{j}_{k=-\infty}\mathcal{P}_{k} H(X_i)(\lambda)\mathcal{P}_{k}H(X_j)(-\lambda) \, d\lambda\\
B_{n,2}
&=\frac{1}{\pi n(n-1)}\sum_{1\leq j<i\leq n, \, i-j>m}\int_{\R} \widehat{K}(\lambda h_n)\sum^{}_{\substack{k\leq i,\, \ell\le j,\, k\neq \ell\\ i-k\leq m_0,\, j-\ell\leq m_0}}\mathcal{P}_{k}H(X_i)(\lambda)\mathcal{P}_{\ell}H(X_j)(-\lambda)\, d\lambda\\
B_{n,3}
&=\frac{1}{\pi n(n-1)}\sum_{1\leq j<i\leq n, \, i-j>m}\int_{\R} \widehat{K}(\lambda h_n)\sum^{}_{\substack{ k\leq i,\, \ell\le j,\, k\neq \ell\\ i-k> m_0\, \text{or}\, j-\ell>m_0}}\mathcal{P}_{k}H(X_i)(\lambda)\mathcal{P}_{\ell}H(X_j)(-\lambda)\, d\lambda.
\end{align*}

Using similar arguments as in {\bf Step 1.}, we can show that 
\begin{align} \label{Bn1e}
\mathbb{E}|B_{n,1}|\leq c_{18}\, n^{1-(\alpha-\eta)\beta}.
\end{align} 
Note that 
\begin{align} \label{Bn2}
&\mathbb{E}|B_{n,2}|^2 \nonumber\\
&\leq \frac{c_{19}}{n^4}\sum_{\substack{1\leq j_{\theta}<i_{\theta}\leq n, \, i_{\theta}-j_{\theta}>m, \, k_{\theta}\neq \ell_{\theta}\\ 0\leq i_{\theta}-k_{\theta}\leq m_0,\, 0\leq j_{\theta}-\ell_{\theta}\leq m_0, \, \theta=1,2}}\int_{\R^2} \big| \widehat{K}(\lambda_1 h_n)\big|\big|\widehat{K}(\lambda_2 h_n)\big| \nonumber\\
&\quad\times \left| \mathbb{E}\big[\mathcal{P}_{k_1}H(X_{i_1})(\lambda_1) \mathcal{P}_{\ell_1}H(X_{j_1})(-\lambda_1)  \mathcal{P}_{k_2}H(X_{i_2})(\lambda_2)  \mathcal{P}_{\ell_2}H(X_{j_2})(-\lambda_2) \big]\right| \, d\lambda_1d\lambda_2.
\end{align}

In the sequel, we will estimate the expectation 
\begin{align} \label{expectation}
\mathbb{E}\big[\mathcal{P}_{k_1}H(X_{i_1})(\lambda_1) \mathcal{P}_{\ell_1}H(X_{j_1})(-\lambda_1)  \mathcal{P}_{k_2}H(X_{i_2})(\lambda_2)  \mathcal{P}_{\ell_2}H(X_{j_2})(-\lambda_2) \big]
\end{align}
and specify the choice of $m$. Assume that $m$ is larger than $4m_0$. Then there are four possibilities for the orderings of $i_1, j_1, i_2,j_2$: 
\[
(1)\; i_1\geq i_2>j_1\geq j_2,  \quad (2)\; i_1\geq i_2>j_2\geq j_1,\quad (3)\; i_2\geq i_1>j_1\geq j_2, \quad (4)\;  i_2\geq i_1>j_2\geq j_1.
\]
By symmetry, it suffices to consider the first two cases. In the first case $i_1\geq i_2>j_1\geq j_2$, the expectation $\eref{expectation}$ is equal to zero if $k_1\neq k_2$.  When $k_1=k_2$,  $0\leq i_1-k_1\leq m_0$ and $0\leq i_2-k_2\leq m_0$ imply $0\leq i_1-i_2\leq m_0$.  If $m>m_0+m_2$, then there is a factor  
\[
\prod^{i_2-m_0-1}_{q=i_2-m_0-m_2}\phi_{\varepsilon}(a_{i_1-q}\lambda_1+a_{i_2-q}\lambda_2)
\]
in the expectation $\eref{expectation}$. By Assumption ({\bf A1}), we can choose $m_2\in\mathbb{N}$ independent of $i_2$ such that 
\begin{align*}
\prod^{i_2-m_0-1}_{q=i_2-m_0-m_2}\left|\phi_{\varepsilon}(a_{i_1-q}\lambda_1+a_{i_2-q}\lambda_2)\right|
&\leq \prod^{m_0+m_2}_{p=m_0+1}\frac{c_{20}}{1+|a_{i_1-i_2+p}\lambda_1+a_{p}\lambda_2|^{\delta}}\\
&\leq \sum^{m_0+m_2}_{p=m_0+1}\frac{c_{21}}{1+|a_{i_1-i_2+p}\lambda_1+a_{p}\lambda_2|^4},
\end{align*}
where in the last inequality we used $\prod\limits^{m_2}_{k=1} x_k\leq \sum\limits^{m_2}_{k=1} x^{m_2}_k$ for any $x_1\geq 0,\cdots, x_{m_2}\geq 0$.

Moreover, if  $|\ell_1-\ell_2|>2m_0+m_3+m_4$, then there is another factor 
\[
\prod^{j_1-m_0-m_3}_{q=j_1-m_0-m_3-m_4}\phi_{\varepsilon}(a_{i_1-q}\lambda_1+a_{i_2-q}\lambda_2-a_{j_1-q}\lambda_1)=\prod^{m_0+m_3+m_4}_{p=m_0+m_3}\phi_{\varepsilon}((a_{i_1-j_1+p}-a_p)\lambda_1+a_{i_2-j_1+p}\lambda_2)
\]
in the expectation $\eref{expectation}$. 

By Assumption ({\bf A1}), $\lim\limits_{n\to\infty}n^{\beta}a_n=c_0$, $|\ell_1-\ell_2|>2m_0+m_3+m_4$ and $0\leq j_{\theta}-\ell_{\theta}\leq m_0$ for $\theta=1,2$, we can choose $m_3, m_4\in\mathbb{N}$ independent of $i_1,i_2,j_1$ such that 
\[
\prod^{j_1-m_0-m_3}_{q=j_1-m_0-m_3-m_4}\left|\phi_{\varepsilon}(a_{i_1-q}\lambda_1+a_{i_2-q}\lambda_2-a_{j_1-q}\lambda_1)\right|\leq \sum^{m_0+m_3+m_4}_{q=m_0+m_3}\frac{c_{22}}{1+|(a_{i_1-i_2+i_2-j_1+q}-a_q)\lambda_1+a_{i_2-j_1+q}\lambda_2|^4}.
\]
Note that 
\[
i_2-j_1= i_1-j_1-(i_1-i_2)>m-m_0.
\]
So we can choose $m\in\mathbb{N}$ large enough such that 
\begin{align} \label{lambda1}
|a_{i_1-i_2+i_2-j_1+q}-a_q|\geq \frac{1}{2}|a_q|
\end{align}
and 
\begin{align} \label{lambda2}
\left|\det\left(\left(\begin{array}{ll}
a_{i_1-i_2+p}  &   a_{p}  \\
a_{i_1-i_2+i_2-j_1+q}-a_q   &   a_{i_2-j_1+q}
\end{array}\right)\right)\right|\geq \frac{1}{4}|a_pa_q|>0
\end{align}
for all $0\leq i_1-i_2\leq m_0$, $m_0+1\leq p\leq m_0+m_2$, $m_0+m_3\leq q\leq m_0+m_3+m_4$.

So for $m$ large enough, in the first case $i_1\geq i_2>j_1\geq j_2$, the right hand side of $\eref{Bn2}$ is less than a constant multiple of 
\[
\text{III}_1+\text{III}_2,
\] 
where
\begin{align*}
\text{III}_1&=\frac{1}{n^4}\sum_{\substack{1\leq j_{\theta}<i_{\theta}\leq n,\, i_{\theta}-j_{\theta}>m,\, k_{\theta}\neq \ell_{\theta}\\ 0\leq i_{\theta}-k_{\theta}\leq m_0,\, 0\leq j_{\theta}-\ell_{\theta}\leq m_0,\, \theta=1,2}}\int_{\R^2} \big| \widehat{K}(\lambda_1 h_n)\big|\big|\widehat{K}(\lambda_2 h_n)\big|\\
&\qquad\qquad\times 1_{\{k_1=k_2, \, |\ell_1-\ell_2|\leq 2m_0+m_3+m_4\}}\sum^{m_0+m_2}_{p=m_0+1}\frac{1}{1+|a_{i_1-i_2+p}\lambda_1+a_{p}\lambda_2|^4}\, d\lambda_1d\lambda_2
\end{align*}
and
\begin{align*}
\text{III}_2&=\frac{1}{n^4}\sum_{\substack{1\leq j_{\theta}<i_{\theta}\leq n,\, i_{\theta}-j_{\theta}>m, \, k_{\theta}\neq \ell_{\theta}\\ 0\leq i_{\theta}-k_{\theta}\leq m_0,\, 0\leq j_{\theta}-\ell_{\theta}\leq m_0,\, \theta=1,2}}\int_{\R^2}1_{\{k_1=k_2,\, |\ell_1-\ell_2|>2m_0+m_3+m_4\}}\sum^{m_0+m_2}_{p=m_0+1}\frac{1}{1+|a_{i_1-i_2+p}\lambda_1+a_{p}\lambda_2|^4}\\
&\qquad\qquad\times \sum^{m_0+m_3+m_4}_{q=m_0+m_3}\frac{1}{1+|(a_{i_1-i_2+i_2-j_1+q}-a_q)\lambda_1+a_{i_2-j_1+q}\lambda_2|^4}\\
&\qquad\qquad\qquad\times \bigg(1_{\{\ell_1<\ell_2\}} |\lambda_1a_{i_1-\ell_2}+\lambda_2a_{i_2-\ell_2}|^{\frac{\alpha-\eta}{2}}|\lambda_1a_{i_1-\ell_1}+\lambda_2a_{i_2-\ell_1}-\lambda_2a_{j_2-\ell_1}|^{\frac{\alpha-\eta}{2}}\\
&\qquad\qquad\qquad\qquad+1_{\{\ell_1>\ell_2\}} |\lambda_1a_{i_1-\ell_1}+\lambda_2a_{i_2-\ell_1}|^{\frac{\alpha-\eta}{2}}|\lambda_1a_{i_1-\ell_2}+\lambda_2a_{i_2-\ell_2}-\lambda_2a_{j_1-\ell_2}|^{\frac{\alpha-\eta}{2}}\bigg) d\lambda_1\, d\lambda_2.
\end{align*}
Clearly, 
\begin{align*}
\text{III}_1
&\leq \frac{1}{n^4}\sum_{\substack{1\leq j_{\theta}<i_{\theta}\leq n,\, i_{\theta}-j_{\theta}>m,\, k_{\theta}\neq \ell_{\theta}\\ 0\leq i_{\theta}-k_{\theta}\leq m_0,\, 0\leq j_{\theta}-\ell_{\theta}\leq m_0,\, \theta=1,2}}\int_{\R^2} \Big(\big| \widehat{K}(\lambda_1 h_n)\big|^2+\big|\widehat{K}(\lambda_2 h_n)\big|^2\Big)1_{\{k_1=k_2, \, |\ell_1-\ell_2|\leq 2m_0+m_3+m_4\}}\\
&\qquad\qquad\qquad\times \sum^{m_0+m_2}_{p=m_0+1}\frac{1}{1+|a_{i_1-i_2+p}\lambda_1+a_{p}\lambda_2|^4}\, d\lambda_1d\lambda_2\\
&\leq \frac{c_{23}}{n^4}\Big(\int_{\R} \big| \widehat{K}(\lambda_1 h_n)\big|^2\, d\lambda_1+\int_{\R} \big| \widehat{K}(\lambda_2 h_n)\big|^2\, d\lambda_2\Big)\sum_{\substack{1\leq j_{\theta}<i_{\theta}\leq n,\, i_{\theta}-j_{\theta}>m,\, k_{\theta}\neq \ell_{\theta}\\ 0\leq i_{\theta}-k_{\theta}\leq m_0,\, 0\leq j_{\theta}-\ell_{\theta}\leq m_0,\, \theta=1,2}}1_{\{k_1=k_2, \, |\ell_1-\ell_2|\leq 2m_0+m_3+m_4\}}\\
&\leq \frac{c_{24}}{n^2h_n},
\end{align*}
where we used Plancherel theorem for the kernel function $K$ in the last inequality.

Moreover,
\begin{align*}
\text{III}_2
&\leq \frac{2}{n^4}\sum_{\substack{1\leq j_{\theta}<i_{\theta}\leq n,\, i_{\theta}-j_{\theta}>m, \, k_{\theta}\neq \ell_{\theta}\\ 0\leq i_{\theta}-k_{\theta}\leq m_0,\, 0\leq j_{\theta}-\ell_{\theta}\leq m_0,\, \theta=1,2}} 1_{\{k_1=k_2,\, |\ell_1-\ell_2|>2m_0+m_3+m_4\}}\int_{\R^2}(|\lambda_1|^{\alpha-\eta}+|\lambda_2|^{\alpha-\eta})\\
&\quad\times \sum^{m_0+m_2}_{p=m_0+1}\frac{1}{1+|a_{i_1-i_2+p}\lambda_1+a_{p}\lambda_2|^4}\sum^{m_0+m_3+m_4}_{q=m_0+m_3}\frac{1}{1+|(a_{i_1-i_2+i_2-j_1+q}-a_q)\lambda_1+a_{i_2-j_1+q}\lambda_2|^4}\\
&\qquad\times \bigg(1_{\{\ell_1<\ell_2\}} (|a_{i_1-\ell_2}|^{\frac{\alpha-\eta}{2}}+|a_{i_2-\ell_2}|^{\frac{\alpha-\eta}{2}})(|a_{i_1-\ell_1}|^{\frac{\alpha-\eta}{2}}+|a_{i_2-\ell_1}|^{\frac{\alpha-\eta}{2}}+|a_{j_2-\ell_1}|^{\frac{\alpha-\eta}{2}})\\
&\qquad\quad+1_{\{\ell_1>\ell_2\}} (|a_{i_1-\ell_1}|^{\frac{\alpha-\eta}{2}}+|a_{i_2-\ell_1}|^{\frac{\alpha-\eta}{2}})(|a_{i_1-\ell_2}|^{\frac{\alpha-\eta}{2}}+|a_{i_2-\ell_2}|^{\frac{\alpha-\eta}{2}}+|a_{j_1-\ell_2}|^{\frac{\alpha-\eta}{2}})\bigg) d\lambda_1d\lambda_2\\
&\leq \frac{c_{25}}{n^4}\sum_{\substack{1\leq j_{\theta}<i_{\theta}\leq n,\, i_{\theta}-j_{\theta}>m, \, k_{\theta}\neq \ell_{\theta}\\ 0\leq i_{\theta}-k_{\theta}\leq m_0,\, 0\leq j_{\theta}-\ell_{\theta}\leq m_0,\, \theta=1,2}} 1_{\{k_1=k_2,\, |\ell_1-\ell_2|>2m_0+m_3+m_4\}}\\
&\qquad\times \bigg(1_{\{\ell_1<\ell_2\}} (|a_{i_1-\ell_2}|^{\frac{\alpha-\eta}{2}}+|a_{i_2-\ell_2}|^{\frac{\alpha-\eta}{2}})(|a_{i_1-\ell_1}|^{\frac{\alpha-\eta}{2}}+|a_{i_2-\ell_1}|^{\frac{\alpha-\eta}{2}}+|a_{j_2-\ell_1}|^{\frac{\alpha-\eta}{2}})\\
&\qquad\quad+1_{\{\ell_1>\ell_2\}} (|a_{i_1-\ell_1}|^{\frac{\alpha-\eta}{2}}+|a_{i_2-\ell_1}|^{\frac{\alpha-\eta}{2}})(|a_{i_1-\ell_2}|^{\frac{\alpha-\eta}{2}}+|a_{i_2-\ell_2}|^{\frac{\alpha-\eta}{2}}+|a_{j_1-\ell_2}|^{\frac{\alpha-\eta}{2}})\bigg)\\
&\leq \frac{c_{26}}{n^{1+(\alpha-\eta)\beta}},
\end{align*}
where in the second inequality we used $\eref{lambda1}$ and $\eref{lambda2}$ to make proper change of variables to get the finiteness of the integral with respect to $\lambda_1$ and $\lambda_2$.

Therefore, for $m$ large enough, in the first case $i_1\geq i_2>j_1\geq j_2$, the right hand side of $\eref{Bn2}$ is less than a constant multiple of $\frac{1}{n^2h_n}+\frac{1}{n^{1+(\alpha-\eta)\beta}}$. Similarly, for $m$ large enough, in the second case $i_1\geq i_2>j_2\geq j_1$, we can also show that the right hand side of $\eref{Bn2}$ is less than a constant multiple of $\frac{1}{n^2h_n}+\frac{1}{n^{1+(\alpha-\eta)\beta}}$. Hence,
\begin{align} \label{Bn2e}
\mathbb{E}|B_{n,2}|^2\leq c_{27} \left(\frac{1}{n^2h_n}+\frac{1}{n^{1+(\alpha-\eta)\beta}}\right).
\end{align}

Now we estimate $\mathbb{E}|B_{n,3}|^2$. Note that 
\begin{align}\label{Bn3}
&\mathbb{E}|B_{n,3}|^2 \nonumber\\
&\leq \frac{c_{28}}{n^4}\int_{\R^2}  \left(\sum_{\substack{1\leq j_{\theta}<i_{\theta}\leq n,\, i_{\theta}-j_{\theta}>m,\, k_{\theta}>\ell_{\theta}\\ i_{\theta}-k_{\theta}>m_0\, \text{or}\, j_{\theta}-\ell_{\theta}>m_0,\, \theta=1,2}}+\sum_{\substack{1\leq j_{\theta}<i_{\theta}\leq n,\, i_{\theta}-j_{\theta}>m,\, k_{\theta}<\ell_{\theta}\\ i_{\theta}-k_{\theta}>m_0\, \text{or}\, j_{\theta}-\ell_{\theta}>m_0,\, \theta=1,2}}\right)\big| \widehat{K}(\lambda_1 h_n)\big|\big|\widehat{K}(\lambda_2 h_n)\big| \nonumber\\
&\qquad\times \Big| \mathbb{E}\big[\mathcal{P}_{k_1}H(X_{i_1})(\lambda_1) \mathcal{P}_{\ell_1}H(X_{j_1})(-\lambda_1)  \mathcal{P}_{k_2}H(X_{i_2})(\lambda_2)  \mathcal{P}_{\ell_2}H(X_{j_2})(-\lambda_2) \big]\Big| \, d\lambda_1d\lambda_2.
\end{align}

Recall the choice of $m_0$ in the proof of Lemma $\ref{lem2}$. Then, by using Cauchy-Schwartz inequality and Lemma \ref{lem1}, we can show that the absolute value of the expectation in $\eref{Bn3}$ is less than a constant multiple of 
\begin{align*}
&\frac{1_{\{k_1>\ell_1, k_2>\ell_2, k_1=k_2\}}}{(1+|\lambda_1|^6)(1+|\lambda_2|^6)}|\lambda_1a_{i_1-k_1}|^{\frac{\alpha-\eta}{2}}|\lambda_2a_{i_2-k_2}|^{\frac{\alpha-\eta}{2}}\Big(1_{\{\ell_1=\ell_2\}}|\lambda_1a_{j_1-\ell_1}|^{\frac{\alpha-\eta}{2}}|\lambda_2a_{j_2-\ell_2}|^{\frac{\alpha-\eta}{2}}\\
&\qquad+1_{\{\ell_1<\ell_2\}} |\lambda_2a_{j_2-\ell_2}|^{\frac{\alpha-\eta}{2}} |\lambda_1a_{i_1-\ell_2}+\lambda_2a_{i_2-\ell_2}|^{\frac{\alpha-\eta}{2}}|\lambda_1a_{j_1-\ell_1}|^{\frac{\alpha-\eta}{2}} |\lambda_1a_{i_1-\ell_1}+\lambda_2a_{i_2-\ell_1}-\lambda_2a_{j_2-\ell_1}|^{\frac{\alpha-\eta}{2}}\\
&\qquad+1_{\{\ell_1>\ell_2\}} |\lambda_1a_{j_1-\ell_1}|^{\frac{\alpha-\eta}{2}}|\lambda_1a_{i_1-\ell_1}+\lambda_2a_{i_2-\ell_1}|^{\frac{\alpha-\eta}{2}}|\lambda_2a_{j_2-\ell_2}|^{\frac{\alpha-\eta}{2}}|\lambda_1a_{i_1-\ell_2}+\lambda_2a_{i_2-\ell_2}-\lambda_2a_{j_1-\ell_2}|^{\frac{\alpha-\eta}{2}}\Big)\\
&\quad+\frac{1_{\{\ell_1>k_1, \ell_2>k_2, \ell_1=\ell_2\}}}{(1+|\lambda_1|^6)(1+|\lambda_2|^6)}|\lambda_1a_{i_1-\ell_1}|^{\frac{\alpha-\eta}{2}}|\lambda_2a_{i_2-\ell_2}|^{\frac{\alpha-\eta}{2}}\Big(1_{\{k_1=k_2\}}|\lambda_1a_{j_1-k_1}|^{\frac{\alpha-\eta}{2}}|\lambda_2a_{j_2-k_2}|^{\frac{\alpha-\eta}{2}}\\
&\qquad+1_{\{k_1<k_2\}} |\lambda_2a_{j_2-k_2}|^{\frac{\alpha-\eta}{2}} |\lambda_1a_{i_1-k_2}+\lambda_2a_{i_2-k_2}|^{\frac{\alpha-\eta}{2}}|\lambda_1a_{j_1-k_1}|^{\frac{\alpha-\eta}{2}} |\lambda_1a_{i_1-k_1}+\lambda_2a_{i_2-k_1}-\lambda_2a_{j_2-k_1}|^{\frac{\alpha-\eta}{2}}\\
&\qquad+1_{\{k_1>k_2\}} |\lambda_1a_{j_1-k_1}|^{\frac{\alpha-\eta}{2}}|\lambda_1a_{i_1-k_1}+\lambda_2a_{i_2-k_1}|^{\frac{\alpha-\eta}{2}}|\lambda_2a_{j_2-k_2}|^{\frac{\alpha-\eta}{2}}|\lambda_1a_{i_1-k_2}+\lambda_2a_{i_2-k_2}-\lambda_2a_{j_1-k_2}|^{\frac{\alpha-\eta}{2}}\Big)\\
&\leq \frac{1_{\{k_1>\ell_1, k_2>\ell_2, k_1=k_2\}}}{(1+|\lambda_1|^2)(1+|\lambda_2|^2)}|a_{i_1-k_1}|^{\frac{\alpha-\eta}{2}}|a_{i_2-k_2}|^{\frac{\alpha-\eta}{2}}\Big(1_{\{\ell_1=\ell_2\}}|a_{j_1-\ell_1}|^{\frac{\alpha-\eta}{2}}|a_{j_2-\ell_2}|^{\frac{\alpha-\eta}{2}}\\
&\qquad+1_{\{\ell_1\neq \ell_2\}} |a_{j_2-\ell_2}|^{\frac{\alpha-\eta}{2}} (|a_{i_1-\ell_2}|^{\frac{\alpha-\eta}{2}}+|a_{i_2-\ell_2}|^{\frac{\alpha-\eta}{2}})|a_{j_1-\ell_1}|^{\frac{\alpha-\eta}{2}} (|a_{i_1-\ell_1}|^{\frac{\alpha-\eta}{2}}+|a_{i_2-\ell_1}|^{\frac{\alpha-\eta}{2}}+|a_{j_2-\ell_1}|^{\frac{\alpha-\eta}{2}})\\
&\quad+\frac{1_{\{\ell_1>k_1, \ell_2>k_2, \ell_1=\ell_2\}}}{(1+|\lambda_1|^2)(1+|\lambda_2|^2)}|a_{i_1-\ell_1}|^{\frac{\alpha-\eta}{2}}|a_{i_2-\ell_2}|^{\frac{\alpha-\eta}{2}}\Big(1_{\{k_1=k_2\}}|a_{j_1-k_1}|^{\frac{\alpha-\eta}{2}}|a_{j_2-k_2}|^{\frac{\alpha-\eta}{2}}\\
&\qquad+1_{\{k_1\neq k_2\}} |a_{j_2-k_2}|^{\frac{\alpha-\eta}{2}} (|a_{i_1-k_2}|^{\frac{\alpha-\eta}{2}}+|a_{i_2-k_2}|^{\frac{\alpha-\eta}{2}})|a_{j_1-k_1}|^{\frac{\alpha-\eta}{2}} (|a_{i_1-k_1}|^{\frac{\alpha-\eta}{2}}+|a_{i_2-k_1}|^{\frac{\alpha-\eta}{2}}+|a_{j_2-k_1}|^{\frac{\alpha-\eta}{2}})\Big).
\end{align*}
Therefore, after simple calculations, we have
\begin{align} \label{Bn3e}
\mathbb{E}|B_{n,3}|^2\leq c_{29} \, n^{2-2(\alpha-\eta)\beta}.
\end{align}

Combining $\eref{B}$, $\eref{Bn1e}$, $\eref{Bn2e}$ and $\eref{Bn3e}$ gives
\begin{align}  \label{Bne}
\mathbb{E}|B_n-\mathbb{E} B_n|\leq c_{30} \left(n^{1-(\alpha-\eta)\beta}+\frac{1}{\sqrt{n^2h_n}}\right).
\end{align}

\noindent
{\bf Step 4.} We estimate $\E\big[|N_n-\overline{N}_n|^2\big]$ where
\[
\overline{N}_n=\frac{1}{2\pi}\int_{\R} \widehat{K}(0) \big(\phi_n(\lambda)-\phi(\lambda) \big) \phi(-\lambda) \,  d\lambda.
\]
Let 
\[
\widetilde{N}_{n}=\frac{1}{2\pi}\int_{\R} \widehat{K}(\lambda h_n) \big(\phi_n(\lambda)-\phi(\lambda) \big) \phi(-\lambda) \,  d\lambda.
\] 
Recall the definition of $N_n$ in $\eref{decomp}$. $|N_n-\widetilde{N}_{n}|$ is less than a constant multiple of $\frac{1}{n}$. Moreover, by Cauchy-Schwartz inequality and Lemma \ref{lem2},
\begin{align*}
\E\big[|\widetilde{N}_n-\overline{N}_n|^2\big]
&\leq \E\left[ \Big|\frac{1}{2\pi}\int_{\R} \big(\widehat{K}(\lambda h_n)-\widehat{K}(0)\big)\big(\phi_n(\lambda)-\phi(\lambda) \big) \phi(-\lambda) \,  d\lambda\Big|^2\right]\\
&\leq \left(\int_{\R} \big|\widehat{K}(\lambda h_n)-\widehat{K}(0)\big|^2|\phi(\lambda)|\, d\lambda \right) \left(\int_{\mathbb{R}} \mathbb{E}  |\phi_n(\lambda)-\phi(\lambda)|^2|\phi(\lambda)| \,  d\lambda \right)\\
&\leq c_{31} \, h^4_n\, n^{1-(\alpha-\eta)\beta}.
\end{align*}
Hence 
\begin{align}  \label{Nne}
\E\big[|N_n-\overline{N}_n|^2\big]\leq c_{32} \left(\frac{1}{n^2}+h^4_n\, n^{1-(\alpha-\eta)\beta}\right).
\end{align}

\noindent
{\bf Step 5.}  It is easy to see that 
\[
\overline{N}_n=\frac{1}{n} \sum\limits^n_{i=1}\left(f(X_i)-\int_{\R} f^2(x)\, dx\right).
\] 
Finally, combining $\eref{decomp}$, $\eref{Ane}$, $\eref{Bne}$ and $\eref{Nne}$ gives
\[
\E\Big|T_n(h_n)-\E T_n(h_n)-\frac{1}{n}\sum^n_{i=1}Y_i\Big|\leq c_{33}\left(n^{1-(\alpha-\eta)\beta}+\frac{1}{\sqrt{n^3h^2_n}}+\frac{1}{\sqrt{n^2h_n}}+h^2_nn^{\frac{1-(\alpha-\eta)\beta}{2}} \right).
\]
This finishes the proof of Theorem \ref{thm1}.
\qed

Finally, we give the proof of Theorem \ref{thm2}. 

\noindent
{\it Proof of Theorem \ref{thm2}}  According to Theorem \ref{thm1}, we only need to consider the asymptotic behavior of $\overline{N}_n$. In the region $1<\alpha<2$, \eref{asym1} and \eref{asym2} follow from Corollary 2.3 in \cite{Koul and Surgailis-2001} and Theorem 2.2 in \cite{Surgailis-2002}, respectively. In the region $0<\alpha<1$, \eref{asym3} follows from Theorem 2.1 in \cite{Honda-2009b} and the paragraph after it.
\qed

\bigskip

$\begin{array}{cc}
\begin{minipage}[t]{1\textwidth}
{\bf Hui Liu}\\
School of Statistics, East China Normal University, Shanghai 200262, China \\
\texttt{lhui56@163.com}
\end{minipage}
\hfill
\end{array}$

\medskip

$\begin{array}{cc}
\begin{minipage}[t]{1\textwidth}
{\bf Fangjun Xu}\\
Key Laboratory of Advanced Theory and Application in Statistics and Data Science - MOE, School of Statistics, East China Normal University, Shanghai, 200062, China \\
NYU-ECNU Institute of Mathematical Sciences at NYU Shanghai, 3663 Zhongshan Road North, Shanghai, 200062, China\\
\texttt{fjxu@finance.ecnu.edu.cn, fangjunxu@gmail.com}
\end{minipage}
\hfill
\end{array}$

\end{document}